\documentclass[12pt]{article}
\usepackage{graphicx, amsmath,amsfonts,amsthm, euscript}
\usepackage{verbatim}
\usepackage{color}
\usepackage{ulem}

\newtheorem{theorem}{Theorem}[section]

\newtheorem{lemma}{Lemma}[section]
\newtheorem{corollary}{Corollary}[section]

\usepackage{graphicx, amsmath}
\usepackage{color}
\definecolor{Blue}{rgb}{0.3,0.3,0.9}

\begin{document}
\normalsize
\title{ Convergence of  quantum random walks with decoherence  }

\maketitle

\vskip 6pt
Shimao Fan$^*$, Zhiyong Feng$^*$, Sheng Xiong$^{**}$ and  Wei-Shih Yang$^*$
\vskip 3pt
$^*$Department of Mathematics\\
        Temple University,
        Philadelphia, PA 19122

$^{**}$Department of Mathematics and Sciences\\
        Edward Waters College,
        Jacksonville, FL 32209

\vskip 6pt
Email: shimao.fan@temple.edu, zhiyong.feng@temple.edu \\
 sheng.xiong@ewc.edu, yang@temple.edu

\vskip 12pt

KEY WORDS:  Quantum walk, Decoherence, Limiting distribution

\vskip12pt

PACS numbers: 05.30.-d, 03.67.Lx, 05.40.-a

\begin{abstract}
In this paper, we study the discrete-time quantum random walks on a line subject to decoherence.  The convergence of  the rescaled position probability distribution $p(x,t)$  depends mainly on  the spectrum of the superoperator $\mathcal{L}_{kk}$.   We show that if 1 is an eigenvalue of the superoperator with multiplicity one and  there is no other eigenvalue   whose modulus equals to 1, then   $\hat {P}(\frac{\nu} {\sqrt t},t)$ converges to a convex combination of normal distributions. In terms of position space, the rescaled probability mass function $p_t (x, t) \equiv p(\sqrt t x, t)$, $ x \in Z/\sqrt t$, converges in distribution to a continuous convex combination of normal distributions. We give an necessary and sufficient condition for a $U(2)$  decoherent quantum walk that satisfies the eigenvalue conditions.
  We also give a complete description of the behavior of quantum walks whose eigenvalues do not satisfy these assumptions. Specific examples such as  the Hadamard walk, walks under real and complex rotations are illustrated. For the $O(2)$ quantum random walks, an explicit formula is provided for the scaling limit  of $p(x,t)$ and their moments. We also obtain exact critical exponents for their moments at the critical point and show universality classes with respect to these critical exponents.

\end{abstract}

\section{Introduction}

\setcounter{equation}{0}

In recent years quantum walks (QWs), as the quantum analog of the classical random walks (CRWs), have attracted great attention from  mathematicians, computer scientists, physicists and engineers. Two forms of  QWs, continuous-time  QWs (CTQW) \cite{FG} and discrete-time  QW (DTQW)\cite{ADZ, MS, ABN, NV,BCA}, are widely studied. In this work, we restrict our discussion to DTQW. In this case, an extra ``coin" degree of freedom is introduced into the system.   Unlike the classical random
walk,  where   the direction of the particle moves  is determined by the outcome of a ``coin flip",  for quantum random walks  both the ``flip" of the coin and the conditional motion of the particle   given by unitary transformations are needed.

In 2003, Brun, Carteret and Ambainis \cite{BCA} discussed the decoherent Hadamard walk on 1-dimensional integer lattice $Z$, and found the expressions for the first and second moments of the position and showed that in the long time limit the variance grows linearly with time with the diffusive character. However, they did not provide an exact expression for the higher order moments nor the limiting distribution, due to the complicated  forms of the  eigenvalues of the superoperator $\mathcal{L}_{kk}$ and  the difficulty to evaluate the position probability analytically.
In 2004, Grimmett, Janson and Scudo \cite{GJS2004} obtained the scaling limit of quantum random walks without decoherence. Their scaling factor is $\frac{1}{t}$.
Recently in \cite{AVWW2011}, the scaling limit of  a quantum random walk with a Markov controlled coin process converges either with scaling factor $\frac{1}{t}$ or $\frac{1}{\sqrt t}$, depending on eigenvalue conditions of the walk operator.

 In this paper, we consider the model with decoherence operators given by \cite{BCA}. We overcome the difficulties in there to obtain the scaling limit of decoherent quantum random walks, by analyzing  $\hat {P}(\nu, t)$, the characteristic function (Fourier transformation) of the position probability distribution $p(x,t)$.   Here $p(x,t)$, $x \in Z$, denotes the probability  that the  particle is found at position $x$ at time $t$. It turns out that  the convergence of $\hat {P}(\frac{\nu}{\sqrt t},t)$  depends on the spectrum of the superoperator $\mathcal{L}_{kk}$.  We show that if 1 is an eigenvalue of the superoperator with multiplicity one and  there is no other eigenvalue   whose modulus equals to 1, then   $\hat {P}(\frac{\nu} {\sqrt t},t)$ converges to a convex combination of normal distributions. In terms of position space, the rescaled probability mass function $p_t (x, t) \equiv p(\sqrt t x, t)$, $ x \in Z/\sqrt t$, converges in distribution to a continuous convex combination of normal distributions. We give an necessary and sufficient condition for a $U(2)$  decoherent quantum walk that satisfies the eigenvalue conditions. For the $O(2)$ quantum walks such as  the Hadamard walk, and walks under real or complex rotations are discussed. An explicit limiting distribution formula is provided for these walks.

Our  article is organized as follows. In Section 2, we present basic concepts of quantum random walks.   In Section 3, we present our main result about the limit of $p_t(x,t)$.  In Section 4, we give examples that illustrate our results. We obtain the scaling limit and the exact critical exponents for their moments at the critical point. We show that
the decoherent quantum random walks, with coin space unitary transformation  $ U\in O(2)$, $\theta\neq \frac{n\pi}{2}$, $n= 0,1,2,3$, belong to the same universality class with respected to the critical exponents of all moments at their critical points.

\section{The unitary walk on the line and decoherence}
  Consider a general quantum random walk on the 1-dimensional integer lattices $Z$.  To be consistent,  we adapt analogous definitions and notations as those outlined in  \cite{BCA}.  We denote the state space by a Hilbert space $H_{p}\otimes H_2$, where  $H_{p}$ denotes the position space and  $H_{2}$ denotes the coin space. The basis of the position space are $|x>$, where $x\in Z$ and,  the basis of the coin space are $|R>$ and $|L>$. We will assume that the walk starts at the origin. The shift operators in $H_p$ are defined as follows
  \begin{equation}
  S^+|x>=|x+1>,
  \end{equation}
  \begin{equation}
 S^-|x>=|x-1>,
  \end{equation}
  where $S^-$ and $S^+$ are unitary shift operators on the particle position.  Let $P_R, P_L$ be two orthogonal projections on the coin space $H_2$ spanned by $|R>$ or $|L>$, where $P_R+P_L=I$. Let $U$ be a unitary transformation on $H_2$, that acts as  the "flipping" of the coin.   Then  the evolution operator of the quantum random walk is given by
    \begin{equation}
  E=(S^+\otimes P_R+S^-\otimes P_L)(I\otimes U).
  \end{equation}
The eigenvectors $|k>$ of $S^-$, $S^+$ are
\begin{equation}
|k>=\sum_x e^{ikx} |x>,\; k\in [0, 2\pi],
\end{equation}
with eigenvalues
\[
S^+|k>=e^{-ik}|k>,
\]
\begin{equation}
S^-|k>=e^{ik}|k>.
\end{equation}
Therefore, in $|k>$ basis, the evolution operator is
\begin{equation}
E(|k>\otimes |\Phi>)=|k>\otimes(e^{-ik}P_R+e^{ik}P_L)U|\Phi>\equiv|k>\otimes U_k|\Phi>,
\end{equation}
where $U_k=(e^{-ik}P_R+e^{ik}P_L)U$ is also a unitary operator.

  The decoherence on the coin space is defined as follows. Suppose before each unitary transformation acting on the coin, a measurement is performed on the coin. This measurement  is given by a set of operators $\{A_n\}$ on $H_2$ which satisfy
\begin{equation}
\sum_n A^*_nA_n=I.
\end{equation}
Through out this paper, we also assume that the measurement is unital, i.e., it satisfies
\begin{equation}
\sum_n A_nA^*_n=I.
\end{equation}

After the measurement, a density operator $\chi$ on $H_2$ is transformed by
\begin{equation}
\chi\rightarrow\chi^{\prime}=\sum_n A_n\chi A^*_n.
\end{equation}
The general density operator of quantum random walk is then given by
\begin{equation}
\rho=\int\frac{dk}{2\pi}\int\frac{dk^{\prime}}{2\pi}|k><k^{\prime}|\otimes \chi_{kk^{\prime}},
\end{equation}
where $\chi_{kk^{\prime}}\in L(H_2)$, and $L(H_2)$ is a vector space of linear operators on $H_2$.
Then after one step of the evolution under coin space decoherence, the density operator  becomes
\begin{equation}
\rho^{\prime}=\int\frac{dk}{2\pi}\int\frac{dk^{\prime}}{2\pi}|k><k^{\prime}|\otimes \sum_nU_kA_n\chi_{kk^{\prime}}A^*_nU^*_{k^\prime}
\end{equation}

Suppose the quantum walk starts at the state $|0>\otimes |\Phi_0>$,
then the initial state is given by the density operator
\begin{equation}
\rho_0=\int\frac{dk}{2\pi}\int\frac{dk^{\prime}}{2\pi}|k><k^{\prime}|\otimes |\Phi_0><\Phi_0|.
\end{equation}
After  $t$ steps, the state evolves to
\begin{equation}
\rho_t=\int\frac{dk}{2\pi}\int\frac{dk^{\prime}}{2\pi}|k><k^{\prime}|\otimes\sum_{n_1,\ldots,n_t}U_kA_{n_t}\cdots U_kA_{n_1} |\Phi_0><\Phi_0| A^*_{n_1}U^*_{k^{\prime}}\cdots A^*_{n_t}U^*_{k^{\prime}}.\label{eq:rhot}
\end{equation}
  If we define the superoperator  $\mathcal{L}_{kk^{\prime}}$ to be an operator which  maps $L(H_2)$ to $L(H_2)$:
\begin{equation}
\mathcal{L}_{kk^{\prime}}B\equiv\sum_{n}U_kA_{n} B A^*_{n}U^*_{k^{\prime}}, \; \forall B\in L(H_2),
 \label{eq:Lkkprime}
\end{equation}
then
\begin{equation}
\rho_t=\int\frac{dk}{2\pi}\int\frac{dk^{\prime}}{2\pi}|k><k^{\prime}|\otimes \mathcal{L}^t_{kk^{\prime}}|\Phi_0><\Phi_0|.
\end{equation}
The probability of reaching a point $x$ at time $t$ is
\begin{eqnarray}
p(x,t)&=&Tr\{(|x><x|\otimes I)\rho_t\}\nonumber\\
&=&\frac{1}{(2\pi)^2}\int dk\int dk^{\prime}<x|k><k^{\prime}|x>Tr\{ \mathcal{L}^t_{kk^{\prime}}|\Phi_0><\Phi_0|\}\nonumber\\
&=&\frac{1}{(2\pi)^2}\int dk\int dk^{\prime} e^{ix(k-k^{\prime})}Tr\{ \mathcal{L}^t_{kk^{\prime}}|\Phi_0><\Phi_0|\}. \label{eq:prob t}
\end{eqnarray}

\section{The limiting distributions of  quantum walks with decoherence coin}
\setcounter{equation}{0}
Let
\begin{equation}
\hat{P}(\nu,t)\equiv<e^{i\nu x}>_t=\sum_xe^{i\nu x}p(x,t)
\end{equation}
be the characteristic function of $p(x,t)$. By the property of  the   $\delta$ function
\begin{equation}
\frac{1}{2\pi}\sum_x x^m e^{-ix(k-k^{\prime})}=(-i)^m\delta^{(m)}(k-k^{\prime}),
\end{equation}
and (\ref{eq:prob t}), we have
\begin{eqnarray}
<e^{i\nu x}>_t &=& \sum_x e^{i\nu x}p(x,t)\nonumber\\
&=&\sum_x e^{i\nu x}\int \frac{dk}{2\pi}\int \frac{dk^{\prime}}{2\pi} e^{ix(k-k^{\prime})}Tr\{ \mathcal{L}^t_{kk^{\prime}}|\Phi_0><\Phi_0|\}\nonumber\\
&=&\int \frac{dk}{2\pi}\int \frac{dk^{\prime}}{2\pi}\sum_x e^{ -ix(k^{\prime}-k-\nu)}Tr\{ \mathcal{L}^t_{kk^{\prime}}|\Phi_0><\Phi_0|\}\nonumber\\
&=&\int \frac{dk}{2\pi}\int \frac{dk^{\prime}}{2\pi}2\pi \delta( k^{\prime}-k-\nu)Tr\{ \mathcal{L}^t_{kk^{\prime}}|\Phi_0><\Phi_0|\}\nonumber\\
&=&\frac{1}{2\pi}\int dk Tr\{ \mathcal{L}^t_{k, k+\nu}|\Phi_0><\Phi_0|\}.
\end{eqnarray}

 For any initial state $\hat{O} \in L(H_2)$, the generating function of $<e^{i\nu x}>_t$ is given by
\begin{eqnarray}
G(z,\nu)&=&\sum_{t=0}^{\infty}z^t<e^{i\nu x}>_t\\
&=&\frac{1}{2\pi}\int dk\sum_{t=0}^{\infty}z^tTr\{\mathcal{L}^t_{k,k+\nu}\hat{O}\}\\
&=&\frac{1}{2\pi}\int dk Tr\{\frac{1}{I-z\mathcal{L}_{k,k+\nu}}\hat{O}\}, \label{eq:Green function}
\end{eqnarray}
where $|z|<1$ and $\hat{O}\in L(H_2)$.  Note that the generating function is well defined since the spectrum of $\mathcal{L}_{k,k+\nu}$ is less than or equal to $1$ by the following lemma.
\begin{lemma} \label{LKV}
Suppose $U\in U_2(\mathbb{C})$ and the set of operators $\{A_n\}$ is unital. Let $\lambda$ be an eigenvalue of $\mathcal{L}_{k,k+\nu}$, then $|\lambda|\leq 1.$
\end{lemma}
\textbf{Proof:} Define
$
\|\mathcal{L}_{k,k}\|=\sup_{0 \ne \hat{O}\in L(H_2)} \frac{\|\mathcal{L}_{k,k}\hat{O}\|}{\|\hat{O}\|}
$
with
$\|\hat{O}\|^2=Tr(\hat{O}^*\hat{O}),$ then we have \\ $\|\mathcal{L}_{k,k}\|\le 1$ (see ref.\cite{BCJZ}).
\begin{eqnarray*}
\|\mathcal{L}_{k,k+\nu}\hat{O}\|^2&=&Tr((\sum_nU_kA_n\hat{O}A_n^*U_{k+\nu}^*)^*(\sum_nU_{k}A_n\hat{O}A_n^*U_{k+\nu}^*))\\
               &=&Tr(U_{k+\nu}(\sum_nA_n\hat{O}A_n^*)^*U_k^*U_{k}(\sum_nA_n\hat{O}A_n^*)U_{k+\nu}^*)\\
               &=&Tr(\sum_nA_n\hat{O}A_n^*)^*(\sum_nA_n\hat{O}A_n^*)) \\
               &=&\|\mathcal{L}_{k,k}\hat{O}\|^2.
\end{eqnarray*}
This shows that $\|\mathcal{L}_{k,k+\nu}\|=\|\mathcal{L}_{k,k}\| \le 1$
 and  therefore $|\lambda|\leq 1$.

It follows from the above lemma  that $\sum_{t=0}^{\infty}(z\mathcal{L}_{k,k+\nu})^t\hat{O}$ converges in $|z|<1$ and thus $I-z\mathcal{L}_{k,k+\nu}$ does not have any pole inside the disk $|z|<1$.

 It is well known that any $\hat{O}\in L(H_2)$ can be written as a linear combination of Pauli matrices:
\begin{equation}
\hat{O}=r_0I+r_1\sigma_1+r_2\sigma_2+r_3\sigma_3,
\end{equation}
where $\sigma_{1,2,3}=\sigma_{x,y,z}$ are usual Pauli matrices. Hence  $\hat{O}$ can be represented by a column vector
\begin{equation}
\hat{O}=
\begin{pmatrix}
r_0\\
r_1\\
r_2\\
r_3
\end{pmatrix}.
\end{equation}

By Lemma \ref{LKV} and  $1-z\mathcal{L}_{k,k+\nu}\in L(L(H_2))$, we have
\begin{equation}
<e^{i\nu x}>_t=\frac{1}{2\pi i}\oint_{|z|=r<1} \frac{G(z,\nu)}{z^{t+1}}dz,
\end{equation}
   for some $0<r<1$. Let  $A$ be the matrix associated with $1-z\mathcal{L}_{k,k+\nu}$, with respect to the Pauli matrices, then
\begin{equation*}\frac{1}{1-z\mathcal{L}_{k,k+\nu}}\hat{O}=
A^{-1}\begin{pmatrix}
r_0\\
r_1\\
r_2\\
r_3
\end{pmatrix}=\frac{1}{\det (A)}\begin{pmatrix}
A_{11} & A_{21} &A_{31} &A_{41}\\
A_{12} & A_{22} &A_{32} &A_{42}\\
A_{13} & A_{23} &A_{33} &A_{43}\\
A_{14} & A_{24} &A_{34} &A_{44}
\end{pmatrix}\begin{pmatrix}
r_0\\
r_1\\
r_2\\
r_3
\end{pmatrix},
\end{equation*}
where $A_{ij}$ is the cofactor of $A$.

Note that $Tr(\sigma_i)=0 $ for $i=1,2,3$, and $Tr(\sigma_0)=2 $.  So when taking the trace in (\ref{eq:Green function}) only the first row action $h(z,\nu)=A_{11}r_0+A_{21}r_1+A_{31}r_2+A_{41}r_4$ remains. Therefore
\begin{eqnarray}G(z,\nu)=\frac{1}{2\pi}\int dk \frac{2h(z,\nu)}{\det A}. \label{eq: formula G}
\end{eqnarray}
 Let  $L=(l_{ij}(\nu))$ be the matrix representation of  $\mathcal{L}_{k,k+\nu}$ in terms of Pauli matrices.  Then we have the following lemma.
\begin{lemma}\label{le:LKK}
 Suppose $U\in U_2(\mathbb{C})$ and $\{A_n\}$ is unital. Then $\mathcal{L}_{k,k+\nu}$ has the following representation
  \[
\begin{pmatrix}
\cos \nu& \times& \times& \times\\
0& \times& \times&\times\\
0& \times& \times& \times\\
i\sin \nu & \times& \times& \times
\end{pmatrix}.
\]

Moreover, if $\nu=0$, then $\mathcal{L}_{k,k}$ has the following representation

 \[
\begin{pmatrix}
1& 0& 0& 0\\
0& \times& \times&\times\\
0& \times& \times& \times\\
0& \times& \times& \times
\end{pmatrix}.
\]
\end{lemma}
\textbf{Proof:} Let $U\in U_2(\mathbb{C})$. Then $|det U|=1$. Let $det U= e^{i\gamma}$. We consider the normalized operator
\begin{eqnarray}
W=e^{-i\frac{\gamma}{2}}U. \label{eq:W}
\end{eqnarray}
Then $W \in SU_2(\mathbb{C})$. By (\ref{eq:Lkkprime}), $\mathcal{L}_{kk^{\prime}} $ is the same for $U$ and $W$. Therefore, without loss of generality, we may assume that $U \in SU_2(\mathbb{C})$ with the following form
\[U=
\begin{pmatrix}
\alpha&-\overline{\beta}\\
\beta&\overline{\alpha}
\end{pmatrix}
\]
where $\alpha,\beta\in\mathbb{C}$ and $|\alpha|^2 + |\beta|^2 = 1$.
Then
\begin{equation}
U_k=
\begin{pmatrix}
e^{-ik}\alpha &-e^{-ik}\overline{\beta}\\
e^{ik}\beta &e^{ik}\overline{\alpha}
\end{pmatrix},
\end{equation}
and $$\mathcal{L}_{k,k+\nu}\sigma_i=l_{1i}(\nu)\sigma_0+l_{2i}(\nu)\sigma_1+l_{3i}(\nu)\sigma_2+l_{4i}(\nu)\sigma_3=\begin{pmatrix}
l_{1i}(\nu)+l_{4i}(\nu)& -il_{2i}(\nu)+l_{3i}(\nu)\\
il_{2i}(\nu)+l_{3i}(\nu)&l_{1i}(\nu)-l_{4i}(\nu)
\end{pmatrix}, $$ where $\; i=0,1,2,3.$ On the other hand,
\begin{equation}
\mathcal{L}_{k,k+\nu}\sigma_0=\sum_{n}U_kA_{n} I A^*_{n}U^*_{k+\nu}=U_k(\sum_{n}A_{n} I A^*_{n})U^*_{k+\nu}=U_k U^*_{k+\nu}=\begin{pmatrix}
e^{i\nu}&0\\
0& e^{-i\nu}
\end{pmatrix}.
\end{equation}
Hence
\begin{eqnarray*}
l_{11}(\nu)+l_{41}(\nu)&=e^{i\nu},\\
l_{11}(\nu)-l_{41}(\nu)&=e^{-i\nu},\\
-il_{21}(\nu)+l_{31}(\nu)&=0,\\
il_{21}(\nu)+l_{31}(\nu)&=0,
\end{eqnarray*}
 which implies that $l_{11}(\nu)=\cos(\nu)$, $l_{21}(\nu)=l_{31}(\nu)=0$ and $l_{41}(\nu)=i\sin (\nu)$.
 In particular, if $\nu=0, $ then  the first column of $\mathcal{L}_{k,k}$ is $(1,0,0,0)^T$.

  Next we finish our proof by showing $l_{1i}(0)=0$, for $i=2,3,4$. Suppose
\begin{equation*}
\sum_nU_kA_n\sigma_{i-1}A_n^*U_k^*=
\begin{pmatrix}
\tau_1^i& \tau_2^i\\
\tau_3^i &\tau_4^i
\end{pmatrix},
\end{equation*}
then
\begin{eqnarray*}
l_{1i}(0)+l_{4i}(0)=\tau_1^i(0),\\
l_{1i}(0)-l_{4i}(0)=\tau_4^i(0),
\end{eqnarray*}
and hence
\begin{equation}
l_{1i}(0)=\frac{1}{2}(\tau_1^i(0)+\tau_4^i(0))=\frac{1}{2}Tr(\mathcal{L}_{k,k}\sigma_{i-1})=\frac{1}{2}Tr(\sigma_{i-1})=0,
\end{equation}
for $i=2,3,4$. The third equality in the above holds because $\mathcal{L}_{k,k}$ preserves the trace.


\begin{theorem}\label{CONV}
Suppose $U\in U_2(\mathbb{C})$,  the set of operators $\{A_n\}$ is unital, $1$ is an eigenvalue of $\mathcal{L}_{k,k}$ with multiplicity one, and  $|\lambda|< 1$ for any other eigenvalue $\lambda$ of  $\mathcal{L}_{k,k}$. Then
\[
\lim_{t\rightarrow\infty}\hat{P}(\frac{\nu}{\sqrt{t}},t)= \frac{1}{2\pi}\int_0^{2\pi}  e^{-\frac{1}{2}z^{\prime\prime}_0(0)\nu^2} dk, \; \forall \nu \in [0,2\pi],
\]
where $z_0(\nu)$ is the root of $det (1-z\mathcal{L}_{k,k+\nu})=0$ such that $z_0(0)=1$.
\end{theorem}

By the Cramer-Levy Theorem (see e.g. Theorem 6.3.2 \cite{CKL}), the above theorem implies that under the conditions, the distribution of the scaling limit is a continuous convex combination of normal distributions with variance $z^{\prime\prime}_0(0)$. In other words, if we define  the rescaled probability mass function on $\frac{Z}{\sqrt t} $ by
\begin{eqnarray}
p_t(x, t)\equiv p(\sqrt t x, t), x \in \frac{Z}{\sqrt t},
\end{eqnarray}
then $p_t$ converges in distribution to the continuous convex combination of normal distributions whose density function is given by
$$
F(x)=\frac{1}{2 \pi}\int_0^{2\pi}\frac{1}{\sqrt {2\pi z^{\prime\prime}_0(0)} }e^{-\frac{1}{2z^{\prime\prime}_0(0)}x^2}dk, x \in R.
$$

\textbf{Proof:} By (\ref{eq:W}), $U$ can be normalized to an $SU_2$ operator $W$. Note that by (\ref{eq:Lkkprime})and (\ref{eq:prob t}),  $\mathcal{L}_{k k+\nu} $ and $p(x,t)$ are the same for $U$ and $W$. Therefore, without loss of generality, we may assume that $U \in SU_2(\mathbb{C})$.  By (\ref{eq: formula G}) and Cauchy's integral formula,
\begin{eqnarray*}
<e^{i\nu x}>_t&=&\frac{1}{2\pi i}\oint_{|z|=r<1} \frac{G(z,\nu)}{z^{t+1}}dz=\frac{1}{2\pi}\int dk\frac{1}{2\pi i}\oint_{|z|=r<1}\frac{2h(z,\nu)} {z^{t+1} det A }dz.
\end{eqnarray*}
Let $z_0 (\nu),z_1(\nu),z_2(\nu),z_3(\nu)$ be four roots of $\det A=0$. Then $\frac{1}{z_0(\nu)},\frac{1}{z_1(\nu)},\frac{1}{z_2(\nu)},\frac{1}{z_3(\nu)}$ are four eigenvalues of $\mathcal{L}_{k,k+\nu}$. By the assumptions that $1$ is an eigenvalue of $\mathcal{L}_{k,k}$ with multiplicity one, and  $|\lambda|< 1$ for any other eigenvalue $\lambda$ of  $\mathcal{L}_{k,k}$, we may make the following ordering $1=|z_0(0)|<  |z_1(0)|\leq|z_2(0)|\leq|z_3(0)|$. Let $l\,(z,\nu)=\frac{2h(z,\nu)} {z^{t+1} det A }$. By continuity of $z_i(\nu), i=0,1,2,3$, in both variables $\nu $ and $k$, there exist a constant $R>1$,  and a small neighborhood
$V$ of $\nu=0$ such that $|z_0(\nu)|<R<|z_i(\nu)|$  for any $\nu \in \bar {V}$ and $k \in [0, 2\pi]$, for all  $i=1,2,3$. By Cauchy's Residue Theorem,
\begin{eqnarray*}
\frac{1}{2 \pi i}\oint_{|z|=R} \frac{2h(z,\nu)} {z^{t+1}det A    }dz= Res(l,z=0)+Res(l,z=z_0(\nu)).
\end{eqnarray*}
 Note that there exists $t_0$ such that $\frac{\nu}{\sqrt{t}}\in \bar{V}$ for all $\nu \in (0, 2\pi)$  if   $t \ge t_0$. Let
 $g(z,\nu)=\det A$. By compactness of $\{(z, \nu, k); |z|=R, \nu \in \bar {V}, k \in [0, 2\pi]\}$, there exists a constant $C$ such that
$|\frac{2h(z,\frac{\nu}{\sqrt{t}})}{g(z,\frac{\nu}{\sqrt{t}})}|\le C$  on $\{(z, \nu); |z|=R, \nu \in [0, 2\pi ],k \in [0, 2\pi] \}$, for all $t \ge t_0$. Therefore,
\begin{eqnarray*}
|\lim_{t\rightarrow\infty}\frac{1}{2\pi i}\oint_{|z|=R}\frac{1}{z^{t+1}}\frac{2h(z,\frac{\nu}{\sqrt{t}})} {g(z,\frac{\nu}{\sqrt{t}})}dz|\leq \lim_{t\rightarrow\infty}\frac{1}{R^{t+1}}C=0.
\end{eqnarray*}
 Hence
$$
\lim_{t\rightarrow\infty}Res(l,z=0)=-\lim_{t\rightarrow\infty}Res(l,z=z_0(\frac{\nu}{\sqrt{t}}))
.$$
For any fixed $\nu$, we have
\[
Res(l,z=z_0(\frac{\nu}{\sqrt{t}}))=\frac{2h(z_0(\frac{\nu}{\sqrt{t}}),\frac{\nu}{\sqrt{t}})}{z_0^{t+1}(\frac{\nu}{\sqrt{t}})\frac{\partial g}{\partial z}(z_0(\frac{\nu}{\sqrt{t}}),\frac{\nu}{\sqrt{t}})}.
\]
Let $t\rightarrow \infty$, we have
\[
\lim_{t\rightarrow\infty}Res(l,z=z_0(\frac{\nu}{\sqrt{t}}))=\frac{2h(1,0)}{\frac{\partial g}{\partial z}(1,0)}\lim_{t\rightarrow\infty}z_0\left(\frac{\nu}{\sqrt{t}}\right)^{-t-1}
\]

We claim that
\begin{eqnarray}
\frac{2h(1,0)}{\frac{\partial g}{\partial z}(1,0)}=-1;\;\; z^{\prime}(0)=0. \label{eq:first derivative of z}
\end{eqnarray}
By the Dominated Convergence Theorem, we have
\begin{eqnarray}
\lim_{t\rightarrow\infty}\hat{P}(\frac{\nu}{\sqrt{t}},t)=\frac{1}{2\pi}\int_0^{2\pi}  e^{-\frac{1}{2}z^{\prime\prime}_0(0)\nu^2} dk.
\end{eqnarray}

We will finish our proof by  proving the claim (\ref{eq:first derivative of z}). Let $M$ denote the 3 by 3 submatrix of $L=(l_{ij})$:
\[
M(\nu)=\begin{pmatrix}
l_{22}&l_{23}&l_{24}\\
l_{32}&l_{33}&l_{34}\\
l_{42}&l_{43}&l_{44}
\end{pmatrix}.
\]

Then the matrix \[A|_{\nu=0}=
\begin{pmatrix}
1-z&\textbf{0} \\
\textbf{0}&I_3-zM(0)
\end{pmatrix}.
\]
The cofactor $A_{21}=A_{31}=A_{41}=0$
and $A_{11}=\det(I_3-zM(0))$.
By Lemma \ref{le:LKK},
 $\frac{1}{z_i(0)}$ are eigenvalues of $\mathcal{L}_{k,k}$ for $i=0,1,2,3$. Therefore $\frac{1}{z_1(0)},\frac{1}{z_2(0)},\frac{1}{z_3(0)}$ are eigenvalues of $M(0)$.
   Hence
\begin{eqnarray}\det(I_3-zM(0))=(1-\frac{z}{z_1(0)})(1-\frac{z}{z_2(0)})(1-\frac{z}{z_3(0)})\\
=-\frac{1}{z_1(0)z_2(0)z_3(0)}(z-z_1(0))(z-z_2(0))(z-z_3(0)).
\end{eqnarray}
On the other hand,
\begin{eqnarray}
\frac{\partial g}{\partial z}(1,0)=\frac{1}{z_1(0)z_2(0)z_3(0)}(1-z_1(0))(1-z_2(0))(1-z_3(0)).\label{eq:partial gz}
\end{eqnarray}
Hence
\begin{eqnarray*}
\frac{2h(1,0)}{\frac{\partial g}{\partial \nu}(1,0)}=2r_0\frac{-\frac{1}{z_1(0)z_2(0)z_3(0)}(1-z_1(0))(1-z_2(0))(1-z_3(0))}{\frac{1}{z_1(0)z_2(0)z_3(0)}(1-z_1(0))(1-z_2(0))(1-z_3(0))}
=-2r_0=-Tr(\hat{O})=-1,
\end{eqnarray*}
since $\hat{O}$ is a density operator.
Next we will show that $z_0^{\prime}(0)=0$. Since $g(z_0(\nu),\nu)=0$, we have
\begin{eqnarray*}
0=\frac{d g(z_0(\nu),z)}{d\nu}=\frac{\partial g(z,\nu)}{\partial \nu}|_{z=z_0(\nu)}+\frac{\partial g(z,\nu)}{\partial z}|_{z=z_0(\nu)}z_0^{\prime}(\nu).
\end{eqnarray*}
When  $\nu=0$,  $z_0(0)=1$, then the above equation becomes
\begin{equation}
0=\frac{\partial g(1,\nu)}{\partial \nu}|_{\nu=0}+\frac{\partial g(z,0)}{\partial z}|_{z=1}z_0^{\prime}(0). \label{eq:implicity derivative}
\end{equation}
Consider the matrix $A$ at $z=1,$ and note   that  $l_{21}(\nu)=l_{31}(\nu)=0$, we have
\begin{eqnarray*}
&&g(1,\nu)=(1-l_{11}(\nu))A_{11}(\nu)+l_{21}(\nu)A_{21}(\nu)
-l_{31}(\nu)A_{31}(\nu)+l_{41}(\nu)A_{41}(\nu)\\
&=&(1-\cos(\nu))A_{11}(\nu) +i\sin {\nu}A_{41}(\nu).
\end{eqnarray*}
By Lemma \ref{le:LKK}, the cofactor  $A_{41}(0)=0$. It follows that
\begin{equation}
\frac{\partial g(1,\nu)}{\partial \nu}|_{\nu =0}=0.\label{eq:partial gnu}
\end{equation}

 By (\ref{eq:partial gz}),  $\frac{\partial g(z_0,\nu)}{\partial z}|_{z=z_0} \neq 0, $ since $z_i(0) \ne  1$, for all $i=1, 2, 3$, by the assumptions of the theorem. By (\ref{eq:implicity derivative}) and (\ref{eq:partial gnu}) we have  $z^{\prime}_0(0)=0$.
\begin{flushright}
\textbf{Q.E.D}
\end{flushright}

\section{Applications}

The assumptions that  1 is the largest eigenvalue  of $\mathcal{L}_{kk}$ with algebraic multiplicity 1 and that there is
no other eigenvalues  whose modulus equals to 1, are crucial in determining the convergence of $ \hat{P}(\frac{\nu}{\sqrt{t}},t)$. In this section, we consider the measurements given by
\[A_0=\sqrt{1-p}I,
\]
\[
A_1=\sqrt{p}|R><R|,
\]
\[
A_2=\sqrt{p}|L><L|.
\]
By analyzing the spectrum of $\mathcal{L}_{kk}$,  we obtain a necessary and sufficient conditions in which the assumptions of Theorem \ref{CONV} are satisfied. Specific examples such as  the Hadamard walk, walks under real and complex rotations are illustrated. For certain  class of convergent quantum walks, explicit formulas are obtained for the limits of the characteristic functions of properly scaled $p(x,t)$. In addition, we will also give a complete description of the behavior of those walks that do not satisfy the conditions.

\begin{lemma}\label{EIG}
Let $\mathcal{L}_{k,k^{\prime}}$ be a superoperator on the Hilbert space $L(\mathbb{C}^2)$, defined by
\[
\mathcal{L}_{k,k^{\prime}}(\hat{O})=\sum_{n=0}^2U_kA_n\hat{O}A_n^*U_{k^{\prime}}^*,
\]
where $U_k$ and $U_{k^{\prime}}$ are  $2\times 2$ unitary matrices and $\hat{O}\in L(\mathbb{C}^2)$.
Then
\[<\mathcal{L}_{k,k^{\prime}}\hat{O},\mathcal{L}_{k,k^{\prime}}\hat{O}>\leq<\hat{O},\hat{O}>.
\]
 In particular, $<\mathcal{L}_{k,k^{\prime}}\hat{O},\mathcal{L}_{k,k^{\prime}}\hat{O}>=<\hat{O},\hat{O}>$ if and only if the decoherence rate $p=0$ or  $\hat{o}_{12}=\hat{o}_{21}=0$.
\end{lemma}

Part of the above lemma has been  obtained by Liu and Petulante in \cite{LP}, but our lemma extends the equality part of their lemma to a wider scope with more applications.

\textbf{Proof}: The inequality follows from Lemma 1 in \cite{LP}. By the proof of Lemma 1 in \cite{LP}, the equality holds  if and only if $(2p-p^2)(\hat{o}_{12}^2+\hat{o}_{21}^2)=0$.  That is, $p=0$ or $\hat{o}_{12}=\hat{o}_{21}=0$.

If  $p\neq 0$ and $\hat{o}_{12}=\hat{o}_{21}=0$, then 1 may not be the unique eigenvalue of $\mathcal{L}_{k,k}$ with largest modulus. In this case, we have the following theorem.  Let $dim(\lambda)$ denote the dimensions of the eigenspace associated with the eigenvalue $\lambda$.

\begin{theorem}\label{SU2} Let $0<p<1$. Let $U\in U(2)$.
Suppose $\lambda$ is an eigenvalue of $\mathcal{L}_{kk}$,
 then we have \\
 (a) $\sum_{|\lambda|=1} dim (\lambda)\le 2.$\\
 (b) 1 is an eigenvalue of  $\mathcal{L}_{k,k}$ and $dim (1)\ge 1.$\\
 (c) If $|\lambda |=1$, then  $\lambda=1 \;\text{or}\; -1$.\\
 (d)  $dim (1) =2$ if and only if  $u_{12}=u_{21}=0$ and $|u_{11}|=|u_{22}|=1$. In this case its multiplicity is 2.\\
(e) There exists eigenvalue $\lambda=-1$ if and only if $u_{11}=u_{22}=0$ and $|u_{12}|=|u_{21}|=1$. In this case its multiplicity is 1. \\

\end{theorem}

\textbf{Proof}: By considering the normalized operator $W$ as in (\ref{eq:W}), and noting that the statements (a)-(e) do not depend on whether it is $U$ or $W$, we mat assume without loss of generality that $U$ is in $SU_2$ with the form
\[U=
\begin{pmatrix}
\alpha&-\overline{\beta}\\
\beta&\overline{\alpha}
\end{pmatrix},
\]
where $\alpha,\beta\in\mathbb{C}$ and $|\alpha|^2 + |\beta|^2 = 1$.

a) Let $\lambda$ be an eigenvalue of $\mathcal{L}_{kk}$ with eigenvector $\hat{O}=\begin{pmatrix}
\hat{o}_{11}&\hat{o}_{12}\\
\hat{o}_{21}&\hat{o}_{22}
\end{pmatrix}.
$By Lemma \ref{le:LKK},
 $\mathcal{L}_{k,k}$ has the form
\begin{eqnarray}
\begin{pmatrix}
1&0&0&0\\
0& l_{22}& l_{23}& l_{24}\\
0& l_{32}&l_{33}& l_{34}\\
0&l_{42}&l_{43}& l_{44}
\end{pmatrix}, \label{eq:LKK form}
\end{eqnarray}
where
\begin{eqnarray}
 l_{24}=-2\cos(2k)Re(\beta\overline{\alpha})+2\sin(2k)Im(\beta\overline{\alpha}),\label{eq:eigen equation 1}\\
l_{34}=-2\cos(2k)Im(\beta\overline{\alpha})-2\sin(2k)Re(\beta\overline{\alpha}),\label{eq:eigen equation 2}\\
l_{42}=2qRe(\alpha \beta),\label{eq:eigen equation 3}\\
l_{43}=2iqIm( \alpha \beta),\label{eq:eigen equation 4}\\
l_{44}=|\alpha|^2-|\beta|^2. \label{eq:eigen equation 5}
\end{eqnarray}
 Since $|\lambda|=1$, by Lemma \ref{EIG},  we have $\hat{o}_{21}=\hat{o}_{12}=0$. This implies that the dimension of  the space spanned by the eigenspace  for all eigenvalues with modulus 1  is at most 2. Moreover, the intersections of eigenspace corresponding to different eigenvalues is $\{0\}$. Therefore a) holds.

b)  By (\ref{eq:LKK form}), 1 is an eigenvalue of  $\mathcal{L}_{k,k}$. Furthermore, $(1,0,0,0)^T$ is one of its eigenvectors. Therefore $dim(1)\ge 1$.

c) Note that
\begin{eqnarray}
\hat{O}=\frac{1}{2}(\hat{o}_{11}+\hat{o}_{22})\sigma_0+\frac{i}{2}(\hat{o}_{12}-\hat{o}_{21})\sigma_1+\frac{1}{2}(\hat{o}_{12}+\hat{o}_{21})\sigma_2+\frac{1}{2}(\hat{o}_{11}-\hat{o}_{22})\sigma_3, \label{eq:pauli rep}
\end{eqnarray}
 so  if $L_{k,k}\hat{O}=\lambda\hat{O}$, with $|\lambda|=1$, then by Lemma \ref{EIG},
\begin{eqnarray}
\frac{1}{2}(\hat{o}_{11}+\hat{o}_{22})=\frac{\lambda}{2}(\hat{o}_{11}+\hat{o}_{22}),\label{eq:eigen equation 6}\\
(2\cos(2k)Re(\beta\overline{\alpha})-2\sin(2k)Im(\beta\overline{\alpha}))\frac{1}{2}(\hat{o}_{11}-\hat{o}_{22})=0,\label{eq:eigen equation 7}\\
(2\cos(2k)Im(\beta\overline{\alpha})+2\sin(2k)Re(\beta\overline{\alpha}))\frac{1}{2}(\hat{o}_{11}-\hat{o}_{22})=0,\label{eq:eigen equation 8}\\
(|\alpha^2|-|\beta^2|)\frac{1}{2}(\hat{o}_{11}-\hat{o}_{22})=\frac{\lambda}{2}(\hat{o}_{11}-\hat{o}_{22}). \label{eq:eigen equation 9}
\end{eqnarray}
(\ref{eq:eigen equation 7}) and (\ref{eq:eigen equation 8}) can be written as the following  matrix form
\begin{eqnarray}
\begin{pmatrix}
\cos 2k & -\sin 2k\\
\sin 2k &\cos 2k
\end{pmatrix} (Re\beta\overline{\alpha}, Im\beta\overline{\alpha})^T (\hat{o}_{11}-\hat{o}_{22})=0. \label{eq:eigen equation 10}
\end{eqnarray}
Note that the matrix on the left hand side of (\ref{eq:eigen equation 10}) has determinant 1, for all $k$, hence invertible. Therefore, if $\hat{o}_{11}-\hat{o}_{22}\neq 0,$ then $|\beta\overline{\alpha}|=0$, or equivalently, $\alpha \beta=0$. Consequently, if $\beta=0$, then $\lambda=1$. If $\alpha=0$, then $\lambda=-1$. On the other hand, if $\hat{o}_{11}-\hat{o}_{22}=0,$ then $\hat{o}_{11}+\hat{o}_{22}\neq 0$ (otherwise $\hat{O}=0$ ). In this case, $\lambda=1$ by (\ref{eq:eigen equation 6}). So $\lambda=1\; \text{or}\; -1$.

d)  We first assume $dim (1) = 2$.  This implies that there exists an eigenvector of the form $(a, 0, 0, b)^T$, with $b \ne 0$. Since $(1,0,0,0)^T$ is an eigenvector,
$(0,0,0,1)^T$ is also an eigenvector.  By (\ref{eq:pauli rep}), $(1,0,0,0)^T$ and $(0,0,0,1)^T$ are two eigenvectors corresponding to $\hat{o}_{11}-\hat{o}_{22}= 0$ and $\hat{o}_{11}+\hat{o}_{22}= 0$, respectively. Hence
$|\alpha|^2-|\beta|^2=1$ by (\ref{eq:eigen equation 9}). Therefore $\beta=0$ and $|\alpha|=1$.

Conversely, if $\beta=0$ and $|\alpha|=1$, then $\mathcal{L}_{k,k}$ has the form
\begin{eqnarray}
\begin{pmatrix}
1&0&0&0\\
0& l_{22}& l_{23}&0\\
0& l_{32}&l_{33}&0\\
0&l_{42}&l_{43}&1
\end{pmatrix}.\label{eq:Lkk form 1 1 case}
\end{eqnarray}
It follows from (\ref{eq:eigen equation 3}) and (\ref{eq:eigen equation 4}) that $l_{42}=l_{43}=0$. Therefore, $dim(1)\ge 2$. So $dim (1)=2$ by part a).  Therefore, by (\ref{eq:Lkk form 1 1 case}) again, the multiplicity of 1 is also two, otherwise $dim(1) > 2$ which contradicts a).

e)  If $-1$ is an eigenvalue, then its dimension must be $1$ by part a). Suppose $(a, 0, 0, 1)^T$ is the associated eigenvector, then $a=0$ by (\ref{eq:eigen equation 6}). Therefore  $|\alpha|^2-|\beta|^2=-1$ by (\ref{eq:eigen equation 9}). This implies  $\alpha=0$ and $|\beta|=1$.

 Conversely, if $\alpha=0$ and $|\beta|=1$, then $\mathcal{L}_{k,k}$ has the form
\begin{eqnarray}
\begin{pmatrix}
1&0&0&0\\
0& l_{22}& l_{23}&0\\
0& l_{32}&l_{33}&0\\
0&0&0&-1
\end{pmatrix}.\label{eq:Lkk form 1 -1 case}
\end{eqnarray}
So  dim(-1)=1 by part a) and b).  Therefore, by (\ref{eq:Lkk form 1 -1 case}) again, the multiplicity of $-1$ is also one, otherwise $dim(-1) \ge 2$ which contradicts a).


\begin{corollary}
 If $ U\in O(2)$, i.e. $ U=
\begin{pmatrix}
\cos \theta & -\sin \theta\\
\sin \theta &\cos \theta
\end{pmatrix}$ or
$\begin{pmatrix}
\cos \theta & \sin \theta\\
\sin \theta &-\cos \theta
\end{pmatrix}$ for some $\theta \in [0, 2\pi]$, then

a) $1$ is an eigenvalue of $\mathcal{L}_{k,k}$ with multiplicity one, and  $|\lambda|< 1$ for any other eigenvalue holds if and only if  $\theta\neq \frac{n\pi}{2}$ where $n= 0,1,2,3$.

b) If  $\theta=0,\pi$, 1 is  an eigenvalue of  $\mathcal{L}_{k,k}$  with multiplicity 2.

c) If $\theta=\frac{\pi}{2},\frac{3\pi}{2}$,
$\mathcal{L}_{k,k}$ has eigenvalues  1 and -1, each has multiplicity 1.
\end{corollary}

\textbf{Proof}:   Note that $\theta=0,\pi$ if and only if  $u_{12}=u_{21}=0$ and $|u_{11}|=|u_{22}|=1$., and  $\theta=\frac{\pi}{2},\frac{3\pi}{2}$ if and only if
$u_{12}=u_{21}=0$ and $|u_{11}|=|u_{22}|=1$. Therefore corollary follows.

We   are now ready to discuss examples according to different values of $\theta$.

In  the Hadamard walk, the evolution operator  is given by
\begin{equation}
U_k=\frac{1}{\sqrt{2}}
\begin{pmatrix}
e^{-ik} &e^{-ik}\\
e^{ik} &-e^{-ik}
\end{pmatrix},
\end{equation}
where $ U=
\begin{pmatrix}
\cos \frac{\pi}{4} & \sin \frac{\pi}{4}\\
\sin \frac{\pi}{4} &-\cos \frac{\pi}{4}
\end{pmatrix}$. By Theorem \ref{SU2}, the associated superoperator
$\mathcal{L}_{k,k+\nu}$ satisfies the assumptions in Theorem \ref{CONV}, and  we have
\[z_0^{\prime}(0)=0; \;\;z_0^{\prime\prime}(0)=\frac{1+q^2+2q\cos 2k}{1-q^2},
\]
where $q=1-p$.
By Theorem \ref{CONV}, we have
\begin{eqnarray}
\lim_{t \to \infty}\hat{P}(\frac{\nu}{\sqrt{t}},t)=\frac{1}{2\pi}\int_0^{2\pi}e^{-\frac{1}{2}\frac{1+q^2+2q\cos 2k}{1-q^2}\nu^2}dk\\
=\frac{1}{2\pi}\int_0^{2\pi}e^{-\frac{1}{2}\frac{1+q^2+2q\cos k}{1-q^2}\nu^2}dk,
\end{eqnarray}
by change of variables and periodicity of cosine.

In general, if $ U\in O(2)$ with $\det U=-1$ and
 $\theta\neq \frac{n\pi}{2}$ where $n= 0,1,2,3$. Then
\begin{equation}
U_k=
\begin{pmatrix}
e^{-ik}\cos(\theta) &e^{-ik}\sin(\theta)\\
e^{ik}\sin(\theta) &-e^{ik}\cos(\theta)
\end{pmatrix},
\end{equation}
and
\[
\mathcal{L}_{k,k+\nu}=
\begin{pmatrix}
\cos(\nu) &(1-p)i\sin(\nu)\sin(2\theta) & 0 &i\sin(\nu)\cos(2\theta)\\
0& -(1-p)\cos(2k+\nu)\cos(2\theta)&(1-p)\sin(2k+\nu)&\cos(2k+\nu)\sin(2\theta)\\
0&-(1-p)\sin(2k+\nu)\cos(2\theta)&-(1-p)\cos(2k+\nu)& \sin(2k+\nu)\sin(2\theta)\\
i\sin(\nu)&(1-p)\cos(\nu)\sin(2\theta)&0&\cos(\nu)\cos(2\theta)
\end{pmatrix}
\]
 which satisfies all the assumptions of Theorem \ref{CONV}.   Direct computation gives  \[z_0^{\prime}(0)=0;\;\;
z_0^{\prime\prime}(0)=\frac{1+2q\cos(2k)+q^2}{1-q^2}\cot^2(\theta).
\]
where $q=1-p$. Hence by Theorem \ref{CONV},
\begin{eqnarray}
\lim_{t \to \infty }\hat{P}(\frac{\nu}{\sqrt{t}},t)=\frac{1}{2\pi}\int_0^{2\pi}e^{-\frac{1}{2}\frac{1+2q\cos(2k)+q^2}{1-q^2}\cot^2(\theta)\nu^2}dk\\
=\frac{1}{2\pi}\int_0^{2\pi}e^{-\frac{1}{2}\frac{1+2q\cos(k)+q^2}{1-q^2}\cot^2(\theta)\nu^2}dk, \label{eq:limit distr so2}
\end{eqnarray}
by change of variables and periodicity of cosine.

 Similar calculations also show that (\ref{eq:limit distr so2}) holds for $ U\in O(2)$ with $\det U=1$ and
 $\theta\neq \frac{n\pi}{2}$, $n= 0,1,2,3$.

The n-th moments $M_n$ of the limiting distribution can be calculated from (\ref{eq:limit distr so2}) by using moment generation functions. Let $$\varphi(\nu)=\frac{1}{2\pi}\int_0^{2\pi}e^{-\frac{1}{2}\frac{1+2q\cos(k)+q^2}{1-q^2}\cot^2(\theta)\nu^2}dk.$$ Then
\begin{eqnarray}
\sum_{n=0}^\infty \frac{1}{n!} i^nM_n\nu ^n=\varphi(\nu) \\
= \sum_{n=0}^\infty \frac{1}{n!}(-\frac{1}{2})^n\cot^{2n}(\theta)\nu^{2n}\frac{1}{(1-q^2)^n}\frac{1}{2\pi}\int_0^{2\pi}(1+2q\cos(k)+q^2)^ndk  \label{eq:m derivative}
\end{eqnarray}
It follows that for $ U\in O(2)$ with
 $\theta\neq \frac{n\pi}{2}$, $n= 0,1,2,3$, we have
\begin{eqnarray}
M_{2n+1}=0, n=0, 1, 2, 3,..., \label{eq:odd moments}
\end{eqnarray}
and for even moments,
\begin{eqnarray}
M_{2n}=\frac{(2n)!}{n!} \cot ^{2n}(\theta) \frac{1}{2^n(1-q^2)^n} \frac{1}{2\pi} \int^{2\pi}_0 (1+2q\cos (2k) +q^2)^n dk \label{eq:even moments}\\
=\frac{(2n)!}{n!} \cot ^{2n}(\theta) \frac{1}{2^n(1-q^2)^n} \frac{1}{2\pi} \int^{2\pi}_0 (e^{i2k}+q)^n(e^{-i2k}+q)^n dk \\
=\frac{(2n)!}{n!} \cot ^{2n}(\theta) \frac{1}{2^n(1-q^2)^n} \frac{1}{2\pi} \int^{2\pi}_0 \sum_{l=0}^n \sum _{l^{\prime}=0}^n {n \choose l} e^{il2k}q^{n-l} {n \choose l^{\prime}}e^{-il^{\prime}2k}q^{n-l^{\prime}} dk \\
=\frac{(2n)!}{n!} \cot ^{2n}(\theta) \frac{1}{2^n(1-q^2)^n}  \sum_{l=0}^n {n \choose l}^2 q^{2(n-l)}, n=0,  1, 2,...  \label{eq:even moments sum}
\end{eqnarray}

Therefore we have
\begin{eqnarray}
M_{2n}=\frac{(2n)!}{n!2^n}  (\frac{\cot ^{2}\theta}{1-q^2})^n  T_n(q), n=0,  1, 2,...  \label{eq:even moments sum final form}
\end{eqnarray}
where $T_n(q)$ is a polynomial of $q$ of order $2n$ given by
\begin{eqnarray}
T_n(q)= \sum_{l=0}^n {n \choose l}^2 q^{2l}.  \label{eq:Tn}
\end{eqnarray}

In particular, for Hadamard walk, the second moment of the limiting distribution is given in terms of
\begin{eqnarray}
T_1(q)= 1+q^2.  \label{eq:T2}
\end{eqnarray}
This result agrees with the results given in \cite{BCA}.

Comparing to the well known $2n$-th moment $N_{2n}$ for the normal distribution with mean $0$ and variance $\sigma ^2=\frac{\cot ^{2}\theta}{1-q^2}$,
\begin{eqnarray}
N_{2n}=\frac{(2n)!}{n!2^n}  (\frac{\cot ^{2}\theta}{1-q^2})^n, n=0,  1, 2,...,  \label{eq:normal even moments }
\end{eqnarray}
we see that the scaling limits of the decoherent quantum random walks are not normally distributed if $q\ne 0$. The deviation from the normal distribution gets larger as the even moments gets larger. However the deviations of the $2n$-th moment are by the same factor $T_n(q)$ for all $\theta\neq \frac{j\pi}{2}$, where $j= 0,1,2,3$.

From (\ref{eq:even moments sum final form}), we also obtain the exact critical exponents for $M_{2n}$ at $p=0$:
\begin{eqnarray}
\gamma_{2n}\equiv \lim_{p \to 0}-\frac{\ln M_{2n}}{\ln p}=n.
\end{eqnarray}
This result shows universality in which the critical exponents do not depend on $\theta$ as long as it converges.  In other words,
the coin-space decoherent quantum random walks, with coin space unitary transformation  $ U\in O(2)$, $\theta\neq \frac{n\pi}{2}$, $n= 0,1,2,3$, belong to the same universality class with respected to the critical exponents of all moments as $p \to 0$.

Now we discuss the limiting behavior of the special cases.

a) $\theta=0,\pi$.
If the initial state is $|0>\otimes|R>$, then $E(|0>\otimes|R>)=|1>\otimes|R>$.
If the initial state is $|0>\otimes|L>$, then $E(|0>\otimes|R>)=-|-1>\otimes|L>$. That is, the walk goes either left or right forever. Hence
\[
\hat{P}(\nu,t)=e^{\pm it\nu}
\]
 depending  on the initial conditions.  In general, if the initial state is  $\phi_0=c_R|R>+c_L|L>$, with $|c_R|^2+|c_L|^2=1,$
then
\[
\hat{P}(\frac{\nu}{t},t)=|c_R|^2e^{i\nu}+|c_L|^2e^{-i\nu}.
\]

b) $\theta=\frac{\pi}{2},\frac{3\pi}{2}$.
If the initial state is $|0>\otimes|R>$, then $E(|0>\otimes|R>)=-|1>\otimes|L>$. If the initial state is $|0>\otimes|L>$, then $E(|0>\otimes|R>)=-|-1>\otimes|R>$. That is, the walk switches back and forth between two positions, which is trivial and
  \[
  \lim_{t\rightarrow \infty}\hat{P}(\frac{\nu}{t^{^\kappa}},t)=1,
 \]
for any $\kappa >0$.

\section{Concluding Remarks}

In this paper we consider coin space decoherent quantum random walks with coin space unitary transformation $U$. We prove that under the eigenvalue conditions, the scaling limit of the probability distribution converges in distribution  to a continuous convex combination of normal distributions. An necessary and sufficient condition is obtained for $U$ to satisfy the eigenvalue conditions. For  $U$ in $O(2)$, an exact form of the limiting distribution is given and the moments of all orders are obtained. For this case, the critical exponents are obtained and we show that all $U$ with rotation angles $\theta\neq \frac{n\pi}{2}$, $n= 0,1,2,3$, belong to the same universality class with respected to the critical exponents of all moments as $p \to 0$. Our analysis is based on the characteristic functions of the position distribution and the analysis of eigenvalues, Theorem \ref{SU2} and its corollary, which plays an important role in the applications of our main convergence theorem, Theorem \ref{SU2}.  We believe that a wider class of  universality should hold for general quantum random walks in general d-dimensional lattices with general rotation $U \in U(n)$. For the future research, it would be very interesting to explore and classify their universality classes with respect to their critical points. On the other hand, we have fixed a measurement in our applications, while our general convergence theorem does not depend on the special form as that in our applications.  An interesting problem would be to understand how the general measurements affect the limiting distributions and, especially their universality classes.




\end{document}